\documentclass{amsart}
 \usepackage{amsmath}
 \usepackage{epsfig} 
 \usepackage{graphics} 

\newtheorem{theorem}{Theorem}
\newtheorem{corollary}{Corollary}

\newtheorem{lemma}{Lemma}

\newtheorem{proposition}{Proposition}
\newtheorem{prop}{Proposition}

\newtheorem{problem 1}{Problem 1}
\newtheorem{problem 2}{Problem 2}
\newtheorem{problem 3}{Problem 3}

\newcommand{\bea}{\begin{eqnarray*}}
\newcommand{\eea}{\end{eqnarray*}}

\title[APPROXIMATION OF PARTIALLY SMOOTH FUNCTIONS]
      {APPROXIMATION OF PARTIALLY SMOOTH FUNCTIONS}

\begin{document}

\begin{abstract}
In this paper we discuss approximation of partially smooth
functions by smooth functions. This problem arises naturally in the
study of laminated currents.
 \end{abstract}

\maketitle

\centerline{\scshape  John Erik Forn\ae ss\footnote{The first author is supported by an NSF grant. Keywords: Approximation, currents, test functions.
2000 AMS classification. Primary: 57R30, Secondary: 32U40}, Yinxia Wang and Erlend Forn\ae ss Wold}
\medskip


\medskip

 \medskip

\section{Introduction}

In this paper we discuss approximation of partially smooth functions by smooth functions.
The motivation of this comes from the study of laminated currents \cite{FS2005}.
For  laminated sets the natural test functions are smooth along leaves and continuous
from leaf to leaf. However, for currents the test functions should be smooth in all directions.
Hence such approximation theorems are very important in the study of laminated currents.

\bigskip

We study here first the case of laminations by curves in $\mathbb R^2$. To prove approximation
we need an extra hypothesis which is automatically satisfied in the complex case of this problem in 
$\mathbb C^2$ because
we then have holomorphic motion. However, this is not satisfied  in general and in a forthcoming paper we will discuss counterexamples.

\bigskip

In the subsequent sections we generalize our result to $\mathbb R^3$. First we deal with the case when
we have a lamination of $\mathbb R^3$ by real surfaces and then we consider the case when we
have a lamination by curves. 

\bigskip

We state and prove all our approximation results locally. We obtain global approximations from these
by using a partition of unity. In a forthcoming paper \cite{FWW} we will discuss the case of holomorphic motion in $\mathbb C^2$ with applications to laminated currents.

\section{ Approximation and smoothing on curves in $\mathbb{R}^2$}

We assume that for every $a \in \mathbb{R}$ we have a smooth ${\mathcal C}^1$ graph
$\Gamma_a$ given by $y=f_a(x),x \in \mathbb{R}.$ We assume that $f_a(0)=a$.
We assume that all graphs are disjoint, and that there is a graph through each point
in $\mathbb{R}^2.$ Moreover we assume that the slope function $F(x,y)=f_a'(x)$ if $y=f_a(x)$ of the graphs is a continuous function on $\mathbb R^2.$
Let $\pi:\mathbb{R}^2 \rightarrow \mathbb{R}$ be defined by
$\pi(x,f_a(x))=a.$ Then the above hypothesis implies that $\pi$ is continuous.

\bigskip

In addition we make the following basic assumption:

\bigskip

$(*)$ \ There is a constant $L$ such that for all $a,a'\in \mathbb{R}$ we have that

$$
|f_{a'}'(x)-f_a'(x)| \leq L\cdot |f_{a'}(x)-f_a(x)|
\log \frac{1}{|f_{a'}(x)-f_a(x)|}
$$

\noindent for all $x \in\mathbb{R}.$

\bigskip

We define a class of partially smooth functions:

\bea
{\mathcal A} & := &  \{\phi \in {\mathcal C}(\mathbb{R}^2); \phi_{|\Gamma_a}\in {\mathcal C}^1(\Gamma_a)\; \forall \; a,\\
& & \Phi(x,y):= d/dx[\phi(x,f_a(x))], y=f_a(x) \in {\mathcal C}(\mathbb{R}^2)\}.
\eea

\begin{theorem}
Let $\phi\in \mathcal A.$ Let $K$ be a compact subset of $\mathbb{R}^2$ and let $\epsilon >0$.
If $(*)$ is satisfied then there exists a function $\psi\in {\mathcal C}^1(\mathbb{R}^2)$ such that
for every point $(x,y)=(x,f_a(x))\in K$:
$$
|\psi(x,y)-\phi(x,y)|<\epsilon, |d/dx [\psi(x,f_a(x))]-d/dx [\phi(x,f_a(x))]|<\epsilon.
$$
\end{theorem}

The theorem follows from the following special case for the function $\phi$ which equals
$a$ on $\Gamma_a.$ 
[This is the function $\pi$ above.]

\begin{proposition}
Let $g\in \mathcal A, g(x,f_a(x))=a$. Let $K$ be a compact subset of $\mathbb{R}^2$ and let $\epsilon >0$. If (*) is satisfied then there exists a function $h\in {\mathcal C}^1(\mathbb{R}^2)$ such that
for every point $(x,y)=(x,f_a(x))\in K$:
$$
|h(x,y)-g(x,y)|<\epsilon, |d/dx [h(x,f_a(x))]|<\epsilon.
$$
\end{proposition}

We first show how the Theorem follows from the proposition. We need a Lemma:

\begin{lemma}
Choose $a,a'\in\mathbb{R}$ such that $0<a'-a<1,$ and put $\delta=a'-a.$
Then
$$
 e^{-L\cdot |x|}\log \frac{1}{\delta}\leq \log \frac{1}{f_{a'}(x)-f_a(x)} \leq e^{L\cdot |x|}\log \frac{1}{\delta},
$$
\noindent for all $x \in \mathbb{R}.$
\end{lemma}

\begin{proof}
Note first that if $\phi\in {\mathcal C}^1 (\mathbb{R})$ is a positive differentiable
function that satisfies $|\phi'(x)| \leq L\phi(x)$ for all $x \in \mathbb{R},$
then $e^{-L|x|}\phi(0)\leq \phi(x)\leq e^{L|x|}\phi(0)$ for all $x \in \mathbb{R}.$
Define $\phi(x):= \log \frac{1}{f_{a'}(x)-f_a(x)}.$ By the basic assumption we have that

$$
|\phi'(x)| =\frac{|f'_{a'}(x)-f'_a(x)|}{f_{a'}(x)-f_a(x)} \leq L\cdot
\log \frac{1}{f_{a'}(x)-f_a(x)} = L \cdot \phi(x),
$$

\noindent and so $\log  \frac{1}{\delta}\cdot e^{-L|x|}\leq \phi(x) \leq \log \frac{1}{\delta}\cdot e^{L|x|}.$

\end{proof}

\begin{proof}
We fix a compact set $K$ in $\mathbb{R}^2.$
Fix a positive integer $J.$ For each integer $j$ set $\phi_j(x)=\phi(x,f_{j/J}(x)).$
Set $\Lambda_j(a)= \cos^2(\frac{\pi J}{2}(a-\frac{j}{J}))$ if $(j-1)/J\leq a \leq (j+1)/J$ and $0$ otherwise.
We let
 $$\psi(x,y):= \sum_j\phi_j(x) \Lambda_j \circ \pi(x,y) .$$

Fix $a, j/J \leq a \leq (j+1)/J.$ We have then

\bea
\left| \psi(x,f_a(x))-\phi(x,f_a(x))\right| & \leq & 
\left|\phi_{j}(x)-\phi(x,f_a(x))\right|\Lambda_{j} (a)\\
& + & 
\left|\phi_{j+1}(x)-\phi(x,f_a(x))\right|\Lambda_{j+1} (a)\\
& \leq & 
\left|\phi(x,f_{j/J}(x))-\phi(x,f_a(x))\right|\\
&  + &  
\left|\phi(x,f_{(j+1)/J}(x)-\phi(x,f_a(x))\right|\\
\eea

Using the Lemma it follows that $\psi(x,f_a(x))$ is arbitrarily close to $\phi(x,f_a(x))$
in sup norm.

\bigskip

We also need to compare $\mathcal C^1$ norms on each graph:

\bea
\left|d/dx( \psi(x,f_a(x)))-d/dx(\phi(x,f_a(x)))\right| & \leq & 
\left|d/dx(\phi_{j}(x)-\phi(x,f_a(x)))\right|\Lambda_{j} (a)\\
& + & 
\left|d/dx(\phi_{j+1}(x)-\phi(x,f_a(x)))\right|\Lambda_{j+1} (a)\\
& \leq & 
\left|d/dx(\phi(x,f_{j/J}(x))-\phi(x,f_a(x)))\right|\\
&  + &  
\left|d/dx(\phi(x,f_{(j+1)/J}(x)-\phi(x,f_a(x)))\right|\\
\eea

Applying the Lemma again, it follows  that $\psi(x,f_a(x))$ and $\phi(x,f_a(x))$ are close
in $\mathcal C^1$ norm.

\bigskip

Hence to prove the Theorem we only need to approximate the function $\psi.$
In order to do this we only need to approximate the functions $\Lambda_j \circ \pi.$
However, this is no problem because $\Lambda_j$ is $\mathcal C^1$ and
$\pi$ is approximated by the function $h$ in the Proposition.

\end{proof}

We proceed to prove Proposition 1.
Now let $b \in \mathbb{R}$ such that $a<b<a'$, and define the following function:

$$
d(x)=\frac{f_b(x)-f_a(x)}{f_{a'}(x)-f_a(x)}.
$$

We have that 

$$
(*) |d'(x)| \leq \frac{|f'_b(x)-f'_a(x)|}{f_{a'}(x)-f_a(x)}+\frac{|f'_{a'}(x)-f'_a(x)|}{f_{a'}(x)-f_a(x)}.
$$

\begin{lemma}
If $|f_b(x)-f_a(x)|\geq \frac{1}{4}|f_{a'}(x)-f_a(x)|$ then
$$
|d'(x)| \leq L \log 4+2L \log \frac{1}{\delta} \cdot e^{L|x|}.
$$
\end{lemma}

\begin{proof}
By $(*)$ and our basic assumption together with the fact that $f_b(x)-f_a(x)<f_{a'}(x)-f_a(x)$ we have that

$$
|d'(x)| \leq L \cdot \log \frac{1}{f_b(x)-f_a(x)}+L \cdot \log \frac{1}{f_{a'}(x)-f_a(x)}.
$$

Since $\frac{1}{f_b(x)-f_a(x)}\leq 4\frac{1}{f_{a'}(x)-f_a(x)}$
we have that $$\log \frac{1}{f_b(x)-f_a(x)}\leq \log 4+ \log \frac{1}{f_{a'}(x)-f_a(x)},$$
\noindent and if we use the Lemma we get 

$$
|d'(x)| \leq L \cdot \log 4+2L  \log \frac{1}{\delta} \cdot e^{L|x|}.
$$

\end{proof}

We next prove the proposition.

\begin{proof}
We construct smooth functions $h_\delta$ such that $h_\delta\rightarrow g$ as $\delta \rightarrow 0.$
For any $\delta>0$ we let $a^\delta(j)=j \cdot \delta$ for $j \in \mathbb{Z}.$ Let $\chi: [0,1]\rightarrow
\mathbb{R}$ be a smooth function such that $\chi(t)=1$ for $0 \leq t \leq \frac{1}{4}$ and
$\chi(t)=0$ for $\frac{3}{4}\leq t \leq 1.$ Let $C$ be a constant such that $|\chi'(t)| \leq C$
for all $t \in [0,1].$ 

We first define $h_\delta$ on the graphs $\Gamma_{a^\delta(j)}$ simply by 
${h_\delta}_{|\Gamma_{a^\delta(j)}}\equiv a^\delta(j).$
Next we define $h_\delta$ between two graphs $\Gamma_{a^\delta(j)}$
and $\Gamma_{a^\delta(j+1)}$:

\bea
h_\delta(x,y) & = & a^\delta(j) \cdot \chi\left(  \frac{y-f_{a^\delta(j)}(x)}{f_{a^\delta(j+1)}(x)-
f_{a^\delta(j)}(x)}  \right)\\
& + & 
a^\delta(j+1) \cdot(1- \chi)\left(  \frac{y-f_{a^\delta(j)}(x)}{f_{a^\delta(j+1)}(x)-
f_{a^\delta(j)}(x)}  \right).
\eea

It is clear that $h_\delta \rightarrow h$ in sup norm on $K.$ We have to check that we have convergence
in ${\mathcal C}^1$ norm on each graph.

Fix a point $(x,y)\in \mathbb{R}^2,$ let $b \in \mathbb{R}$ such that $(x,y)$ lies on the graph
$\Gamma_b.$ We differentiate $h_\delta$ along $\Gamma_b$ at $(x,y).$ We have that $(x,y)$
lies between two graphs $\Gamma_{a^\delta(j)}$ and $\Gamma_{a^\delta(j+1)}.$ Now

$$
|\frac{d}{dx}h_\delta(x,f_b(x))|\leq \delta\cdot C \cdot|
\frac{d}{dx} \left(  \frac{f_b-f_{a^\delta(j)}}{f_{a^\delta(j+1)}-f_{a^\delta(j)}}  \right)(x) |.
$$

If $f_b(x)-f_{a^\delta(j)}(x)< \frac{1}{4} (f_{a^\delta(j+1)}(x)-f_{a^\delta(j)}(x))$, we have that $h$ is constant
near $(x,y)$ so that is fine. Otherwise it follows from Lemma 2 that 

$$
\frac{d}{dx} h_\delta(x,f_b(x)) \leq \delta\cdot C \cdot (L \log 4+2\cdot L \log \frac{1}{\delta}
\cdot e^{L\cdot |x|}).
$$

Since this estimate only depends on $|x|$ and $\delta$ the result follows.

\end{proof}

\section{Approximation and smoothing on surfaces in $\mathbb{R}^3$}

We assume that for every $a \in \mathbb{R}$ we have a smooth ${\mathcal C}^1$ surface $\Gamma_a$
given by $z=f_a(x,y), (x,y)\in \mathbb{R}^2$. We assume that all surfaces are disjoint, and that there is a surface through every point in $\mathbb{R}^3.$ Let $\pi: \mathbb{R}^3 \rightarrow \mathbb{R}$ be defined by $\pi(x,y,f_a(x,y))=a.$ Then the above hypotheses imply that $\pi$ is continuous.

In addition we make the following basic assumption:
There is a constant $L$ such that for all $a,a'\in \mathbb{R}$ we have that

$$
(1) \left|\frac{\partial}{\partial x} f_{a'}(x,y)-\frac{\partial}{\partial x}f_a(x,y)\right| \leq L\cdot |f_{a'}(x,y)-f_a(x,y)|
\log \frac{1}{|f_{a'}(x,y)-f_a(x,y)|}
$$
$$
(2) \left|\frac{\partial}{\partial y} f_{a'}(x,y)-\frac{\partial}{\partial y}f_a(x,y)\right| \leq L\cdot |f_{a'}(x,y)-f_a(x,y)|
\log \frac{1}{|f_{a'}(x,y)-f_a(x,y)|}
$$

\noindent for all $x,y \in\mathbb{R}.$

This implies in particular that $F(x,y,z):= \frac{\partial}{\partial x}f_a(x,y)$ if $z=f_a(x,y)$ is a continuous function. Similarly  $G(x,y,z):= \frac{\partial}{\partial y}f_a(x,y)$ if $z=f_a(x,y)$ is a continuous function.

We define a class of partially smooth functions:

\bea
{\mathcal A} & := & \{\phi \in {\mathcal C}(\mathbb{R}^3); \phi_{|\Gamma_a}\in {\mathcal C}^1(\Gamma_a)\; \forall \; a,\\
& & \Phi(x,y,z):= \frac{\partial}{\partial x}\phi(x,y,f_a(x,y)), z=f_a(x,y) \in {\mathcal C}(\mathbb{R}^3),\\
& & \Psi(x,y,z):= \frac{\partial}{\partial y}\phi(x,y,f_a(x,y)), z=f_a(x,y) \in {\mathcal C}(\mathbb{R}^3)\}.
\eea

\begin{theorem}
Let $\phi\in \mathcal A.$ Let $K$ be a compact subset in $\mathbb{R}^3$ and let $\epsilon >0$.
If (1) and (2) are satisfied, then there exists a function $\psi\in {\mathcal C}^1(\mathbb{R}^3)$ such that
for every point $(x,y,z)=(x,y,f_a(x,y))\in K$:
\bea
|\psi(x,y,z)-\phi(x,y,z)| & < & \epsilon,\\
|\frac{\partial}{\partial x}  [\psi(x,y,f_a(x,y))]-\frac{\partial}{\partial x} [\phi(x,y,f_a(x,y))]| & < & \epsilon,\\
|\frac{\partial}{\partial y}  [\psi(x,y,f_a(x,y))]-\frac{\partial}{\partial y} [\phi(x,y,f_a(x,y))]| & < & \epsilon.
\eea
\end{theorem}

As in the previous section this follows from the following result:

\begin{proposition}
Let $g\in \mathcal A, g(x,y,f_a(x,y))=a$. Let $K$ be a compact subset in $\mathbb{R}^3$ and let $\epsilon >0$. Suppose (1) and (2) are satisfied.
Then there exists a function $h\in {\mathcal C}^1(\mathbb{R}^3)$ such that
for every point $(x,y,z)=(x,y,f_a(x,y))\in K$:
\bea
|h(x,y,z)-g(x,y,z)| & < & \epsilon,\\
|\frac{\partial}{\partial x}  [h(x,y,f_a(x,y))]| & < & \epsilon,\\
|\frac{\partial}{\partial y}  [h(x,y,f_a(x,y))]| & < & \epsilon.
\eea
\end{proposition}

We proceed to prove the Proposition.

\begin{proof}
The proof is the same as in the last section:
For any $\delta>0$ we let $a^\delta(j)=j \cdot \delta$ for $j \in \mathbb{Z}.$ Let $\chi: [0,1]\rightarrow
\mathbb{R}$ be a smooth function such that $\chi(t)=1$ for $0 \leq t \leq \frac{1}{4}$ and
$\chi(t)=0$ for $\frac{3}{4}\leq t \leq 1.$ Let $C$ be a constant such that $|\chi'(t)| \leq C$
for all $t \in [0,1].$ 

We first define $h_\delta$ on the surfaces  $\Gamma_{a^\delta(j)}$ simply by 
${h_\delta}_{|\Gamma_{a^\delta(j)}}\equiv a^\delta(j).$
Next we define $h_\delta$ between two surfaces $\Gamma_{a^\delta(j)}$
and $\Gamma_{a^\delta(j+1)}$:

\bea
h_\delta(x,y,z) & = & a^\delta(j) \cdot \chi\left(  \frac{z-f_{a^\delta(j)}(x,y)}{f_{a^\delta(j+1)}(x,y)-
f_{a^\delta(j)}(x,y)}  \right)\\
& + & 
a^\delta(j+1) \cdot(1- \chi)\left(  \frac{z-f_{a^\delta(j)}(x,y)}{f_{a^\delta(j+1)}(x,y)-
f_{a^\delta(j)}(x,y)}  \right).
\eea

It is clear that $h_\delta \rightarrow h$ in sup norm on $K.$ We have to check that we have convergence
in ${\mathcal C}^1$ norm on each surface.

For each $j$ and $a^\delta(j)<b<a^\delta(j+1)$ define as in the last section

$$
d(x,y)=\frac{f_b(x,y)-f_{a^\delta(j)}(x,y)}{f_{a^\delta(j+1)}(x,y)-f_{a^\delta(j)}(x,y)}.
$$

Suppose that $|f_b(x,y)-f_{a^\delta(j)}(x,y)|\geq \frac{1}{4} |f_{a^\delta(j+1)}(x,y)-f_{a^\delta(j)}(x,y)|.$ As in the proof of Lemma 2 we get that

$$
\left| \frac{\partial}{\partial x} d(x,y)\right| \leq L \log 4+ 2 \cdot L \log \frac{1}{f_{a^\delta(j+1)}(x,y)-f_{a^\delta(j)}(x,y)}
$$

\noindent and the same for $ \frac{\partial}{\partial y} d(x,y).$ Using Lemma 1 along the line through $(0,0)$ and $(x,y)$ and increasing $L$ we obtain

$$
\left| \frac{\partial}{\partial x} d(x,y)\right| \leq L \log 4+ 2 \cdot L \log \frac{1}{\delta}\cdot
e^{L \sqrt{x^2+y^2}}.
$$

This gives that

$$
|\frac{\partial}{\partial x} h(x,y,f_b(x,y))|\leq \delta\cdot C \cdot \left(L \log 4
+2 \cdot L \cdot \log \frac{1}{\delta} \cdot e^{L \sqrt{x^2+y^2}}\right).
$$

The same holds for $\frac{\partial}{\partial y}.$

\end{proof}

\section{Approximation and smoothing on Curves in $\mathbb{R}^3$}

We assume that we have for every $a=(a_1,a_2) \in \mathbb{R}^2$ a 
smooth ${\mathcal C}^1$ graph $\Gamma_a$, $y=(y_1,y_2)=f_a(x),x \in \mathbb{R}$. 
We suppose that $f_a(0)=a.$ Also we assume that the graphs are disjoint. We assume 
there is a graph through each point. So for every $(x,y)\in\mathbb{R}^3$ there is a unique vector $a$ 
so that $y=f_a(x).$ We set $F(x,y):=\frac{df_a(x)}{dx}$ when $y=f_a(x).$ 

Our basic assumption is now that :
$$
(BA)  \ \|F(x,y')-F(x,y)\| \leq L\|y'-y\| \log (1/\|y'-y\|)
$$
\noindent for a fixed constant $L$ and for all $x,y,y'$.  For later use we will assume that
$L\log 2>1$.\

\begin{lemma}
Let $F$ satisfy $(BA)$.  There exist positive constants $C,\delta_{0}\in\mathbb{R}$
such that the following holds: For all  positive $\delta\in\mathbb{R}$  with 
$0<\delta<\delta_{0}$ there is a vector valued function $F_{\delta}(x,y)$  such that

\

(i) \ $\|F_{\delta}-F\|_{\infty}\leq C\delta\log\frac{1}{\delta}$, \

(ii) \ $\|\mathrm{J}_yF_{\delta}\|_{\infty}\leq \frac{C}{2}
\mathrm{log}\frac{1}{\delta}$.
 
\end{lemma}
\begin{proof}
We use the function $\Lambda(t)=[\cos (t)]^2, -\pi/2\leq t \leq \pi/2$ and $0$ otherwise.
Next choose a grid, $y_{m,n}= (m\delta, n \delta).$ Let $F_{mn}(x)$ be a smoothing
of the function $F(x,y_{m,n})$ such that  $|F(x,y_{m,n})-F_{mn}(x)|<\delta.$
Next we define 

$$F_\delta(x,y)=\sum_{m,n} \Lambda(\frac{\pi(y_1-m\delta)}{2\delta})\Lambda(\frac{\pi(y_2-n\delta)}{2\delta})
F_{mn}(x).$$

Suppose that $m \delta \leq y_1\leq (m+1)\delta, n \delta \leq y_2\leq (n+1)\delta$.
Then

\bea
& & \|F_\delta(x,y_1,y_2)-F(x,y_1,y_2)\| \\
& = &  \|\sum_{i=m,m+1,j=n,n+1}\Lambda(\frac{\pi(y_1-i\delta)}{2\delta})\Lambda(\frac{\pi(y_2-j\delta)}{2\delta})
F_{ij}(x)
-F(x,y_1,y_2)\|\\
& = &  \|\sum_{i=m,m+1,j=n,n+1}\Lambda(\frac{\pi(y_1-i\delta)}{2\delta})\Lambda(\frac{\pi(y_2-j\delta)}{2\delta})
\left[F_{ij}(x)
-F(x,y_1,y_2)\right]\|\\
& \leq  &  \sum_{i=m,m+1,j=n,n+1}\|
F_{ij}(x)-F(x,y_1,y_2)\|\\
& \leq  &  \sum_{i=m,m+1,j=n,n+1}\|F_{ij}(x)-F(x,i\delta,j\delta)\|+ \|F(x,i\delta,j\delta)-F(x,y_1,y_2)\|\\
& \leq & 4\delta + 4L (2\delta) \log \frac{1}{2\delta}\leq C\delta\log
\frac{1}{\delta},\\
\eea
if $C$ is chosen appropriately. 

\

We estimate the $y$-derivatives of $F_\delta.$

\bea
& & \|\frac{\partial F_\delta}{\partial y_1}\|\\
& = & \| \partial/\partial y_1  \sum_{i=m,m+1,j=n,n+1} \Lambda(\frac{\pi(y_1-i\delta)}{2\delta})
 \Lambda(\frac{\pi(y_2-j\delta)}{2\delta})
F_{ij}(x)\|\\
& = & \| \partial/\partial y_1  \sum_{i=m,m+1,j=n,n+1} \Lambda(\frac{\pi(y_1-i\delta)}{2\delta})
 \Lambda(\frac{\pi(y_2-j\delta)}{2\delta})
[F_{ij}(x)-F_{mn}(x)]\|\\
& \leq & C\log (1/\delta),\\
\eea
if $C$ is chosen appropriately.  Since we can increase $C$ if we like, the lemma is proved. 
\end{proof}

Let $(x,y) \in \mathbb R^3$. Then $y=f_a(x)$ for some $a \in \mathbb{R}^2.$
We define $\pi: \mathbb{R}^3\rightarrow \mathbb{R}^2$,
$\pi(x,y)=a$ if $y=f_a(x).$

\begin{lemma}
The function $\pi$ is continuous.
\end{lemma}

\begin{proof} 
This is proved by the same method as in Lemma 1.
\end{proof}

We introduce again a class of partially smooth functions.

\bea
{\mathcal A} & := &  \{\phi \in {\mathcal C}(\mathbb{R}^3); \phi_{|\Gamma_a}\in {\mathcal C}^1(\Gamma_a)\; \forall \; a,\\
& & \Phi(x,y):= \frac{d}{dx}[\phi(x,f_a(x))], y=f_a(x) \in {\mathcal C}(\mathbb{R}^3)\}.
\eea
We let $D_{R}$  denote the real polydisk of radius $R$: 
$$D_{R}:=
\{(x,y)\in\mathbb{R}^3;|x|\leq R, |y_{i}|\leq R\}.
$$  

\begin{theorem}
Assume (BA) is satisfied.  There exists
a positive $R\in\mathbb{R}$  such that for all $\phi \in\mathcal{A}$ and all $\epsilon > 0$
there exists a ${\mathcal C}^1$ smooth function $\psi:D_{R}
\rightarrow \mathbb{R}$ so that 
$$
|\psi(x,y)- \phi(x,y)|<\epsilon
$$
for all $(x,y)<D_{R}$.  Moreover,
$$\left|\frac{d}{dx}(\psi(x,f_a(x)))-\frac{d}{dx}
(\phi(x,f_a(x)))\right |<\epsilon
$$ 
for all $x,a$ with $(x,f_a(x))\in D_{R}$.
\end{theorem}
As before this follows from the following Proposition:

\begin{prop}
There exists
a positive $R\in\mathbb{R}$  such that for all $\epsilon>0$
there exists a ${\mathcal C}^1$ smooth function $\psi:D_{R}
\rightarrow \mathbb{R}^2$ so that 
$\|\psi(x,y)- \pi(x,y)\|<\epsilon$
for all $(x,y)\in D_{R}$. Moreover, $\left\|\frac{d}{dx}(\psi(x,f_a(x)))\right\|
<\epsilon$ for all $(x,a)$ with $(x,f_a(x))\in D_{R}$.
\end{prop}

\begin{proof}
For small enough 
$\delta$ let $F_{\delta}$  be the approximating function from Lemma 3.\

For each $a\in \mathbb{R}^2$  let $y=f^\delta_a(x)$ be the unique 
(local) solution to the differential equation $dy/dx=F_\delta(x,y)$ with $f^\delta_a(0)=a$.
There is a positive constant $\mu$ such that the solution exist for $|x|<\mu$.
Our first goal is to show that we have $\mu=\infty$. \

We want to compare the graphs of $f^{\delta}_{a}$  and $f_{a}$, so we define a 
function
$$
\phi^{\delta}_{a}(x):=\|f^{\delta}_{a}(x)-f_{a}(x)\|.
$$
We have that

\bea
\|(f^\delta_a)'(x)-f'_a(x)\| & = & 
\|F_\delta(x,f_a^\delta(x))-F(x,f_a(x))\|\\
& = & 
\|F_\delta(x,f_a^\delta(x))-F(x,f^\delta_a(x))+F(x,f_a^\delta(x))-F(x,f_a(x))\|\\
& \leq & 
\|F_\delta(x,f_a^\delta(x))-F(x,f^\delta_a(x))
\|+\|F(x,f_a^\delta(x))-F(x,f_a(x))\|\\
& \leq & \delta+ L\|f_a^\delta(x)-f_a(x)\|\log \frac{1}{\|f_a^\delta(x)-f_a(x)\|}\\
\eea
Differentiating and using the Cauchy-Schwarz inequality we get that
$$
|\frac{\partial}{\partial x}\phi^{\delta}_{a}(x)|\leq\delta+L\phi^{\delta}
_{a}(x)\log\frac{1}{\phi^{\delta}_{a}(x)}.
$$
This gives that $|\frac{\partial}{\partial x}\log\phi^{\delta}_{a}(x)|
\leq \frac{\delta}{\phi^{\delta}_{a}(x)}+L\log\frac{1}{\phi^{\delta}_{a}(x)} $,
and so $f^{\delta}_{a}(x)$  stays bounded for bounded $x$.   Thus $F_{\delta}$ provides us with  a new
lamination of $\mathbb{R}^3$ and our next goal is to show that we actually approximate the 
lamination defined by $F$ on compact sets.  
\begin{lemma}
Let $\phi$ be a ${\mathcal C}^1$ smooth function and suppose $\delta<\phi(x)<1/2$ on $(0,b)$,
 $\phi(0)=\delta$ and $\phi'\leq \delta+L\phi \log (1/\phi)$ on $(0,b).$
If $L \log 2>1,$ then $\phi(x) \leq \delta^{(e^{-2Lx})}$ on $(0,b)$.
\end{lemma}

\begin{proof}
We have $\phi'\leq \delta+L\phi \log (1/\phi)\leq 2L\phi \log (1/\phi).$
Hence 

\bea
\frac{\phi'}{\phi \log (1/\phi)} & \leq & 2L\\
                                 & \Rightarrow & \\
{[\log (\log (1/\phi))]}' & = & \frac{1}{\log (1/\phi)}\frac{-\phi'}{\phi}\\
& \geq & -2L\\
|[\log (\log (1/\phi))]|^x_0 & \geq & -2Lx, x \geq 0\\
\log (\log (1/\phi(x)))-\log (\log (1/\delta)) & \geq & -2Lx\\
\log \left[\frac{\log (1/\phi(x))}{\log (1/\delta)}\right]
& \geq & -2Lx\\
\frac{\log (1/\phi(x))}{\log (1/\delta)} & \geq & e^{-2Lx}\\
\log (1/\phi(x)) & \geq &  \log (1/\delta) e^{-2Lx}\\
\log \phi(x) & \leq & (\log \delta)e^{-2Lx}\\
& = & \log (\delta^{(e^{-2Lx})})\\
\phi(x) & \leq & \delta^{(e^{-2Lx})}\\
\eea

\end{proof}

\begin{corollary}
Let $0<\tau<1$. Then if $|x|\leq (\log (1/\tau))/(2L)$ and if $\delta$ is small enough
such that $\delta^{\tau}<\frac{1}{2}$, we have that
$$
\phi^{\delta}_{a}(x) \leq \delta^\tau.
$$
\end{corollary}

\begin{proof}
If $\phi^{\delta}_{a}(x)\leq \delta,$ we are done. On an interval where
$\phi^{\delta}_a(x)>\delta$ we use the above lemma to show that
$\phi^{\delta}_a(x)\leq \delta^\tau$ as long as
$e^{-2L|x|}\geq \tau,$ i.e. $-2L|x|\geq \log \tau$ so
$|x|\leq (\log (1/\tau))/(2L).$
\end{proof}

This Corollary says that integral curves of $F$ and $F^\delta$ starting at the
same point when $x=0$ remain close.  \

Let $\Gamma^{\delta}$ denote the set
$$
\Gamma^{\delta}=\{(x,f_{a}^{\delta}(x));|x|\leq 1, \|a\|<1\},
$$
and choose an (initial) $R$ such that $D_{R}\subset\Gamma^{\delta}$ as long as 
$\delta$ is small enough.  Define a projection $\pi_{\delta}:D_{R}
\rightarrow\mathbb{R}^2$ by $\pi_{\delta}(x,f^{\delta}_{a}(x))=a$.   Corollary 1
above shows that $\pi_{\delta}\rightarrow\pi$  as $\delta\rightarrow 0$.  \
We will define $\psi=\pi_{\delta}$ for a small enough $\delta$. \
Let $(x_0,y_0)=(x_0,f_{a}(x_0))\in D_{R} $.
It remains to prove that we have
$$
\left\| \frac{d}{dx}\left[\pi_\delta(x,f_a(x))\right]\right\|
<\epsilon
$$
\noindent at $(x_0,y_0)$ if $\delta$ is small enough.  

\bigskip

By the chain rule, if we write $\pi_\delta= (\pi^1_\delta,\pi^2_\delta)$,

$$
\frac{d}{dx}(\pi^i_\delta(x_{0},f_a(x_{0}))=\frac{\partial \pi^i_\delta}
{\partial x}(x_{0},f_a(x_{0}))+\nabla_y\pi^i_\delta(x_{0},f_a(x_{0}))\cdot 
f'_a(x_{0}).
$$

We write this for short as

$$
\frac{d}{dx}(\pi_\delta(x_{0},f_a(x_{0}))=\frac{\partial \pi_\delta}
{\partial x}(x_{0},f_a(x_{0}))+\nabla_y\pi_\delta(x_{0},f_a(x_{0}))\cdot 
f'_a(x_{0}).
$$

Since $\pi_\delta$ is constant on the graphs $y=f^\delta_{a'}(x)$ we also
have $\pi_\delta(x,f_{a'}^\delta(x))\equiv a'$ so

$$0= \frac{\partial \pi_\delta}{\partial x}(x,f_{a'}^\delta(x))+\nabla_y\pi_\delta(x,f_{a'}^\delta(x))\cdot
(f_{a'}^\delta)'(x).$$

Let $a'$ be so that $y_0=f_{a'}^\delta(x_0)=f_a(x_0).$

Subtracting we have

\bea
\frac{d}{dx} (\pi_\delta(x,f_a(x)))(x_0) & = &
\nabla_y \pi_\delta(x_0,y_0)\cdot
\left[f'_a(x_0)-(f^\delta_{a'})'(x_0)\right]\\
& = & 
\nabla_y \pi_\delta(x_0,y_0)\cdot
\left[F(x_0,y_0)-F_\delta(x_0,y_0)\right]\\
\eea

Hence

$$
\left\|\frac{d}{dx}(\pi_\delta(x,f_a(x)))(x_{0})\right\|
\leq C\delta\log\frac{1}{\delta} \left\| \nabla_y \pi_\delta(x_0,y_0)\right\|.
$$
We proceed to prove that if $|x|<\frac{1}{2C}$  then $\| 
\nabla_y \pi_\delta(x_0,y_0)\|<\frac{2}{\sqrt{\delta}}$. \

To this end we consider the nearby graphs $y=f^\delta_{a'}
(x)$ and $y=f^\delta_{a'+\Delta a}(x).$ Define
$(\Delta f^\delta)(x) := f^\delta_{a'+\Delta a}(x)-f_{a'}^\delta(x).$
Then

\bea
\frac{d}{dx}((\Delta f^\delta)(x)) & = & 
(f^\delta_{a'+\Delta a})'(x)-(f^\delta_{a'})'(x)\\
& = & F_\delta(x,f^\delta_{a'+\Delta a}(x))-F_\delta(x,f^\delta_{a'}(x))\\
\eea

By Lemma 3 we get

\bea
\left\|\frac{d}{dx}((\Delta f^\delta)(x))\right\|
 & = & \left\| F_\delta(x,f^\delta_{a'+\Delta a}(x))-
F_\delta(x,f^\delta_{a'}(x))
\right\|\\
& \leq &
C \log (1/\delta)\|f^{\delta}_{a'+\Delta a}(x)-f_{a'}^\delta(x)\|\\
& = & 
C \log (1/\delta)\|(\Delta f^\delta)(x)\|\\
\\
\eea
Since $|\frac{d}{dx}\|(\Delta f^\delta)(x)\||\leq\|\frac{d}{dx}
(\Delta f^\delta)(x)\|$ we get as in the proof of Lemma 1 in Section 1 that
$$
e^{-C\log\frac{1}{\delta}|x|}\|(\Delta f^\delta)(0)\|
\leq \|(\Delta f^\delta)(x)\|\leq e^{C\log\frac{1}{\delta}|x|}
\|(\Delta f^\delta)(0)\|.
$$
But $\Delta f^\delta(0)=f^\delta_{a'+\Delta a} (0)-f^\delta_{a'}(0)=
(a'+\Delta a)-a'=\Delta a$, and so we get 
$$
\delta^{C|x|}\|\Delta a\|\leq\|(\Delta f^\delta)(x)\|\leq\delta^{-C|x|}\|\Delta a\|.
$$

\bigskip

Hence

$$
\frac{\|\Delta a\|}{\|\Delta f^{\delta}(x) \|}\leq\delta^{-C|x|},
$$

and so

 $$ \left\|\frac{\partial \pi^j_\delta}{\partial y_{i}}\right\|
\leq
\frac{1}{\delta^{C|x|}},
$$
for $j,i=1,2$.  If we choose $|x|$ so small that 
$$
C|x|\leq 1/2
$$ 

then

$$
\|\nabla_y\pi_\delta\|
\leq \frac{2}{\sqrt{\delta}}.
$$

But then

$$
\left\| \frac{d}{dx}(\pi_\delta(x,f_a(x)))\right\|
\leq (C\delta\log\frac{1}{\delta})\cdot\frac{2}{\sqrt{\delta}}
= 2C\sqrt{\delta}\log\frac{1}{\delta}.
$$
After decreasing, $R$,  this completes the proof since $\sqrt{\delta}\log\frac{1}{\delta}
\rightarrow 0$  as $\delta\rightarrow 0$.  

\end{proof}

\bigskip

\noindent John Erik Forn\ae ss\\
Mathematics Department\\
The University of Michigan\\
East Hall, Ann Arbor, MI 48109\\
USA\\
fornaess@umich.edu\\

\noindent Yinxia Wang\\
Department of Mathematics\\
Henan University\\
Kaifeng, 475001\\
China\\
yinxiawang@gmail.com\\

\noindent Erlend Forn\ae ss Wold\\
Mathematisches Institut\\
Universitat Bern\\
Sidlerstr. 5\\
CH-3012 Bern\\
Switzerland\\
erlendfw@student.matnat.uio.no\\

\end{document}